\documentclass[12pt,leqno]{article}

\def\diam{\mathop{\rm diam}}

\def\Lip{\mathop{\rm Lip}}

\newtheorem{theorem}{Theorem}
\newtheorem{lemma}[theorem]{Lemma}
\newtheorem{proposition}[theorem]{Proposition}
\newtheorem{sublemma}[theorem]{Sublemma}
\newtheorem{definition}[theorem]{Definition}
\newtheorem{corollary}[theorem]{Corollary}
\newtheorem{problem}[theorem]{Problem}
\newtheorem{remark}[theorem]{Remark}
\newtheorem{claim}[theorem]{Claim}
\newtheorem{assumptions}[theorem]{Assumptions}
\newtheorem{examples}[theorem]{Examples}
\newtheorem{basicfact}[theorem]{Basic Fact}

\newcommand{\begintheorem}{\addtocounter{equation}{1}\begin{theorem}}
\newcommand{\beginlemma}{\addtocounter{equation}{1}\begin{lemma}}
\newcommand{\beginproposition}{\addtocounter{equation}{1}\begin{proposition}}
\newcommand{\beginsublemma}{\addtocounter{equation}{1}\begin{sublemma}}
\newcommand{\begindefinition}{\addtocounter{equation}{1}\begin{definition}}
\newcommand{\begincorollary}{\addtocounter{equation}{1}\begin{corollary}}
\newcommand{\beginproblem}{\addtocounter{equation}{1}\begin{problem}}
\newcommand{\beginremark}{\addtocounter{equation}{1}\begin{remark}}
\newcommand{\beginclaim}{\addtocounter{equation}{1}\begin{claim}}
\newcommand{\beginassumptions}{\addtocounter{equation}{1}\begin{assumptions}}
\newcommand{\beginexamples}{\addtocounter{equation}{1}\begin{examples}}
\newcommand{\beginbasicfact}{\addtocounter{equation}{1}\begin{basicfact}}

\begin{document}

\title{Some topics concerning harmonic analysis on metric spaces}


\author{Stephen Semmes}

\date{}

\maketitle




\begin{abstract}
In this brief survey we give an introduction to some aspects of
``atoms'' on metric spaces and their connection with linear
operators.  
\end{abstract}



	Let $(M, d(x,y))$ be a metric space.  For the purposes of this
article, the underlying set $M$ will always be assumed to contain at
least $2$ distinct elements.  As usual, the distance function $d(x,y)$
should satisfy $d(x,y) \ge 0$ for all $x, y \in M$, $d(x,y) = 0$ if
and only if $x = y$, $d(x,y) = d(y,x)$ for all $x$ and $y$, and the
``triangle inequality''
\begin{equation}
	d(x,z) \le d(x,y) + d(y,z)
\end{equation}
for all $x$, $y$ and $z$ in $M$.  

	Notice that
\begin{equation}
\label{|d(x,z) - d(y,z)| le d(x,y)}
	|d(x,z) - d(y,z)| \le d(x,y) 
\end{equation}
for all $x$, $y$, and $z$ in $M$, which can be derived from the
triangle inequality.  For each $x$ in $M$ and positive real number
$r$, let us write $B(x,r)$ and $\overline{B}(x,r)$ for the open and
closed balls of radius $r$ in $M$, i.e.,
\begin{equation}
	B(x,r) = \{y \in M : d(y,x) < r\}, \ 
		\overline{B}(x,r) = \{y \in M : d(y,x) \le r\}.
\end{equation}
If $E$ is a nonempty subset of $M$, then $\diam E$ denotes the
\emph{diameter} of $E$, defined by
\begin{equation}
	\diam E = \sup \{d(u,v) : u, v \in E\}.
\end{equation}

	Let $s$ be a positive real number.  We say that $(M, d(x,y))$
is \emph{Ahlfors-regular of dimension $s$} if $M$ is complete as a
metric space, and if there is a positive Borel measure $\mu$ on $M$
such that
\begin{equation}
	C_1^{-1} \, r^s \le \mu(\overline{B}(x,r)) \le C_1 \, r^s
\end{equation}
for some positive real number $C_1$, all $x$ in $M$, and all $r > 0$
such that $r \le \diam M$ if $M$ is bounded.  As a basic example, if
$M$ is $n$-dimensional Euclidean space ${\bf R}^n$ with the standard
metric, and if $\mu$ is Lebesgue measure, then in fact
$\mu(\overline{B}(x,r))$ is equal to a constant times $r^n$, where the
constant is simply the volume of the unit ball.  More exotically, one
can consider (simply-connected) nonabelian nilpotent Lie groups, such
as the Heisenberg groups.  These can be given as Euclidean spaces
topologically, but with very different distance functions that are
compatible with the group structure in place of ordinary vector
addition.  For these spaces one still has natural dilations as on
Euclidean spaces, and Lebesgue measure is compatible with both the
group structure and the dilations, in such a way that the measure of a
ball of radius $r$ is equal to a constant times $r^s$, where $s$ is
now a geometric dimension that is larger than the topological
dimension.

	Fix a metric space $(M, d(x,y))$ and a measure $\mu$ on $M$
satisfying the conditions in the definition of Ahlfors-regularity,
with dimension $s$.  The following fact is sometimes useful: there is
a constant $k_1 \ge 1$ so that if $x$ is an element of $M$ and $r$,
$R$ are positive numbers, with $r \le R$, then the ball
$\overline{B}(x,R)$ can be covered by a collection of at most $k_1 \,
(R/r)^s$ closed balls of radius $r$.  If $M$ is bounded, then we may
as well assume that $r < \diam M$ here, because $M$ is automatically
contained in a single ball with radius $\diam M$.  We may also assume
that $R \le \diam M$, since we could simply replace $R$ with $\diam M$
if $R$ is initially chosen to be larger than that.

	To establish the assertion in the preceding paragraph, let us
begin with a preliminary observation.  Suppose that $A$ is a subset of
$\overline{B}(x,R)$ such that $d(x,y) > r$ for all $x$, $y$ in $A$.
Then the number of elements of $A$ is at most $k_1 \, (R/r)^s$, if we
choose $k_1$ large enough (independently of $x$, $R$, and $r$).  Indeed,
\begin{equation}
	\sum_{a \in A} \mu(\overline{B}(a, r/2))
		= \mu\biggl(\bigcup_{a \in A} \overline{B}(a, r/2) \biggr)
		\le \mu(\overline{B}(x, 3R/2)),
\end{equation}
where the first equality uses the disjointness of the balls
$\overline{B}(a, r/2)$, $a \in A$.  The Ahlfors-regularity property
then applies to give a bound on the number of elements of $A$ of the
form $k_1 \, (R/r)^s$.  Now that we have such a bound, suppose that $A$
is also chosen so that the number of its elements is maximal.  Then
\begin{equation}
	\overline{B}(x,R) \subseteq \bigcup_{a \in A} \overline{B}(a,r).
\end{equation}
In other words, if $z$ is an element of $\overline{B}(x,R)$, then
$d(z,a) \le r$ for some $a$ in $A$, because otherwise we could add $z$
to $A$ to get a set which satisfies the same separation condition as
$A$, but which has $1$ more element.  This yields the original
assertion.

	In particular, closed and bounded subsets of $M$ are compact.
This uses the well-known characterization of compactness in terms of
completeness and total boundedness, where the latter holds for bounded
subsets of $M$ by the result just discussed.

	Let us look at some special families of functions on $M$,
called \emph{atoms} (following \cite{CW2}).  For the sake of
definiteness, we make the convention that a ``ball'' in $M$ means a
closed ball (with some center and radius), if nothing else is
specified.  Suppose that $p$ is a real number and $r$ is an extended
real number such that
\begin{equation}
\label{conditions on p, r}
	0 < p \le 1, \ 1 \le r \le \infty, \ p < r.
\end{equation}
An integrable complex-valued function $a(x)$ on $M$ will be called a
\emph{$(p,r)$-atom} if it satisfies the following three conditions:
first, there is a ball $B$ in $M$ such that the support of $a$ is
contained in $B$, i.e., $a(x) = 0$ when $x \in M \backslash B$;
second,
\begin{equation}
\label{int_M a(x) d mu(x) = 0}
	\int_M a(x) \, d \mu(x) = 0;
\end{equation}
and third, 
\begin{equation}
\label{L^r average of a le mu(B)^{-1/p}}
	\biggl(\frac{1}{\mu(B)} \int_M |a(x)|^r \, d \mu(x) \biggr)^{1/r}
		\le \mu(B)^{-1/p}.
\end{equation}
If $r = \infty$, then (\ref{L^r average of a le mu(B)^{-1/p}})
is interpreted as meaning that the supremum (or essential supremum,
if one prefers) of $a$ is bounded by $\mu(B)^{-1/p}$.

	The size condition (\ref{L^r average of a le mu(B)^{-1/p}})
may seem a bit odd at first.  A basic point is that it implies
\begin{equation}
	\int_M |a(x)|^p \, d \mu(x) \le 1,
\end{equation}
by Jensen's inequality.  The index $r$ reflects a kind of regularity
of the atom, and notice that a $(p, r_1)$-atom is automatically a $(p,
r_2)$-atom when $r_1 \ge r_2$.  There are versions of this going in
the other direction, from $r_2$ to $r_1$, and we shall say more about
this soon.

	If $\alpha$ is a positive real number no greater than $1$,
define $\Lip \alpha$ to be the space of complex-valued functions
$\phi(x)$ on $M$ such that
\begin{equation}
	\sup \biggl\{\frac{|\phi(x) - \phi(y)|}{d(x,y)^\alpha} :
			x, y \in M, \ x \ne y \biggr\} < \infty.
\end{equation}
In this case we define $\|\phi\|_{\Lip \alpha}$ to be this supremum.
Notice that $\|\phi\|_{\Lip \alpha} = 0$ if and only if $\phi$
is constant, and that $\|\cdot\|_{\Lip \alpha}$ is a seminorm,
which means that 
\begin{equation}
	\|\phi + \psi\|_{\Lip \alpha} 
		\le \|\phi\|_{\Lip \alpha} + \|\psi\|_{\Lip \alpha}
\end{equation}
and
\begin{equation}
	\|\lambda \, \phi\|_{\Lip \alpha} = |\lambda| \, \|\phi\|_{\Lip \alpha}
\end{equation}
for all $\phi$, $\psi$ in $\Lip \alpha(M)$ and all complex numbers
$\lambda$.

	As in (\ref{|d(x,z) - d(y,z)| le d(x,y)}), $\phi(x) = d(x,z)$
lies in $\Lip 1$ on $M$ for any fixed $z$ in $M$, and indeed
$\|\phi\|_{\Lip 1} = 1$.  If $f$ lies in the analogue of $\Lip 1$ on
the space ${\bf C}$ of complex numbers (with respect to the usual
Euclidean metric), and if $\psi$ is any function in $\Lip \alpha$ on
$M$, then the composition $f \circ \psi$ also lies in $\Lip \alpha$ on
$M$.  One can use this to show that $\Lip 1$ contains ``plenty'' of
nontrivial functions with bounded support in $M$.  This extends to
$\Lip \alpha$ for $\alpha \in (0,1)$, because any function in $\Lip 1$
with bounded support also lies in $\Lip \alpha$ for all $\alpha \in
(0,1)$.

	Suppose that $a(x)$ is a $(p,r)$-atom on $M$ and that $\phi(x)$
lies in $\Lip \alpha$ on $M$.  Consider the integral
\begin{equation}
	\int_M a(x) \, \phi(x) \, d\mu(x).
\end{equation}
Let $B = \overline{B}(z,t)$ be the ball associated to $a(x)$ as
in the definition of an atom.  The preceding integral can be
written as
\begin{equation}
	\int_{\overline{B}(z,t)} a(x) \, (\phi(x) - \phi(z)) \, d\mu(x),
\end{equation}
using also (\ref{int_M a(x) d mu(x) = 0}).  Thus
\begin{eqnarray}
	\Bigl| \int_M a(x) \, \phi(x) \, d\mu(x) \Bigr|
   & \le & \int_{\overline{B}(z,t)} |a(x)| \, |\phi(x) - \phi(z)| \, d\mu(x)
									\\
 & \le & 
   \mu(\overline{B}(z,t))^{1 - (1/p)} \, t^\alpha \, \|\phi\|_{\Lip \alpha}.
								\nonumber
\end{eqnarray}
Ahlfors-regularity implies that
\begin{equation}
	\Bigl| \int_M a(x) \, \phi(x) \, d\mu(x) \Bigr|
  \le C_1^{1 - (1/p)} \, t^{(1 - (1/p)) s + \alpha} \, \|\phi\|_{\Lip \alpha}.
\end{equation}
In particular, 
\begin{equation}
\label{|int_M a(x) phi(x) dmu(x)| le C_1^{1 - (1/p)} ||phi||_{Lip alpha}}
	\Bigl| \int_M a(x) \, \phi(x) \, d\mu(x) \Bigr|
		  \le C_1^{1 - (1/p)} \, \|\phi\|_{\Lip \alpha}
\end{equation}
when $\alpha = ((1/p) - 1) \, s$.

	If we want to be able to choose $\alpha = ((1/p) - 1) \, s$
and have $\alpha \le 1$, then we are lead to the restriction
\begin{equation}
\label{p ge frac{s}{s+1}}
	p \ge \frac{s}{s+1}.
\end{equation}
Indeed, this condition does come up for some results, even if much of
the theory works without it.  There can also be some funny business at
the endpoint, so that one might wish to assume a strict inequality in
(\ref{p ge frac{s}{s+1}}), or some statements would have to be
modified when equality holds.

	In some situations this type of restriction is not really
necessary, perhaps with some adjustments.  Let us mention two basic
scenarios where this happens.  First, suppose that our metric space
$M$ is something like a self-similar Cantor set, such as the classical
``middle-thirds'' Cantor set.  If we define $\Lip \alpha$ on $M$ in
the same way as before, but allowing $\alpha$ to be larger than $1$,
then there are plenty of $\Lip \alpha$ functions, and, for that
matter, there are plenty of functions which are locally constant.  The
computation giving (\ref{|int_M a(x) phi(x) dmu(x)| le C_1^{1 - (1/p)}
||phi||_{Lip alpha}}) still works when $\alpha > 1$, and this is true
in general.  The point is that this naive extension of $\Lip \alpha$
on a metric space $M$ can be degenerate when $\alpha > 1$, e.g., it
may contain only constant functions.  This is true when $M$ is equal
to ${\bf R}^n$ with the standard metric, for instance.  For if $\alpha
> 1$, then any function in $\Lip \alpha$ has derivative $0$
everywhere.

	On the other hand, if $M = {\bf R}^n$ with the standard
Euclidean metric, then there other ways to define classes of more
smooth functions, through conditions on higher derivatives.  In
connection with this, there is a simple way to strengthen (\ref{int_M
a(x) d mu(x) = 0}), which is to ask that the integral of an atom times
a polynomial of degree at most some number is equal to $0$.  If one
does this, then there are natural extensions of (\ref{|int_M a(x)
phi(x) dmu(x)| le C_1^{1 - (1/p)} ||phi||_{Lip alpha}}) for $\alpha >
1$, obtained by subtracting a polynomial approximation to $\phi(x)$.

	A basic manner in which atoms can be used is to test 
localization properties of linear operators.  Suppose that $T$
is a bounded linear operator on $L^2(M)$, and that $a$ is a
$(p,2)$-atom on $M$.  Consider
\begin{equation}
	T(a)
\end{equation}
(as well as $T^*(a)$, for that matter).  This is well-defined as an
element of $L^2(M)$, since $a$ lies in $L^2(M)$.  If $B =
\overline{B}(z,t)$ is the ball associated to $a$ in the definition of
an atom, then the estimate
\begin{eqnarray}
	\enspace
	\biggl(\frac{1}{\mu(B)} \int_M |T(a)(x)|^2 \, d \mu(x) \biggr)^{1/2}
  & \le & \|T\|_{2,2} \, 
	\biggl(\frac{1}{\mu(B)} \int_M |a(x)|^2 \, d \mu(x) \biggr)^{1/2}  
									\\
	& \le & \|T\|_{2,2} \, \mu(B)^{-1/p}			\nonumber
\end{eqnarray}
provides about as much information about $T(a)$ around $B$, on $2 B =
\overline{B}(z,2t)$, say, as one might reasonably expect to have.
However, in many situations one can expect to have decay of $T(a)$
away from $B$, in such a way that
\begin{equation}
\label{||T(a)||_p le k}
	\|T(a)\|_p \le k 
\end{equation}
for some constant $k$ which does not depend on $a$.

	In this argument it is natural to take $r = 2$, but a basic
result in the theory is that one has some freedom to vary $r$.
Specifically, if $b$ is a $(p,r)$-atom on $M$, then it is possible
to write $b$ as
\begin{equation}
\label{b = sum_i beta_i b_i}
	b = \sum_i \beta_i \, b_i,
\end{equation}
where each $b_i$ is a $(p,\infty)$-atom, each $\beta_i$ is a complex
number, and $\sum_i |\beta_i|^p$ is bounded by a constant that does
not depend on $b$ (but which may depend on $p$ or $r$).  Let us give a
few hints about how one can approach this.  As an initial
approximation, one can try to write $b$ as
\begin{equation}
\label{b = beta' b' + sum_j gamma_j c_j}
	b = \beta' \, b' + \sum_j \gamma_j \, c_j,
\end{equation}
where $b'$ is a $(p,\infty)$-atom, $\beta'$ is a complex number such
that $|\beta'|$ is bounded by a constant that does not depend on $b$,
each $c_j$ is a $(p,r)$-atom, and $\sum_j |\gamma_j|^p \le 1/2$, say.
If one can do this, then one can repeat the process indefinitely to
get a decomposition as in (\ref{b = sum_i beta_i b_i}).  In order to
derive (\ref{b = beta' b' + sum_j gamma_j c_j}), the method of
Calder\'on--Zygmund decompositions can be employed.

	Recall that
\begin{equation}
	\Bigl(\sum_k \tau_k \Bigr)^p \le \sum_k \tau_k^p
\end{equation}
for nonnegative real numbers $\tau_k$ and $0 < p \le 1$.
As a consequence, if $\{f_k\}$ is a family of measurable
functions on $M$ such that 
\begin{equation}
	\int_M |f_k(x)|^p \, d\mu(x) \le 1 \qquad\hbox{for all } k,
\end{equation}
and if $\{\theta_k\}$ is a family of constants, then
\begin{equation}
	\int_M \Bigl|\sum_k \theta_k \, f_k(x) \Bigr|^p \, d\mu(x)
		\le \sum_k |\theta_k|^p.
\end{equation}
Because of this, bounds on $\sum_l |\alpha_l|^p$ are natural when
considering sums of the form $\sum_l \alpha_l \, a_l$, where the
$a_l$'s are $(p,r)$-atoms and the $\alpha_l$'s are constants.

	A fundamental theorem concerning atoms is the following.
Suppose that $T$ is a bounded linear operator on $L^2(M)$ again.  (One
could start as well with a bounded linear operator on some other $L^v$
space, with suitable adjustments.)  Suppose also that there is a
constant $k$ so that (\ref{||T(a)||_p le k}) holds for all
$(p,2)$-atoms, where $0 < p \le 1$, as before, or even simply for all
$(p,\infty)$-atoms.  Then $T$ determines a bounded linear operator
on $L^q$ for $1 < q < 2$.  This indicates how atoms are sufficiently
abundant to be useful.

	The proof of this theorem relies on an argument like the one
in Marcinkeiwicz interpolation.  In the traditional setting, one of
the main ingredients is to take a function $f$ in $L^q$ on $M$, and,
for a given positive real number $\lambda$, write it as $f_1 + f_2$,
where $f_1(x) = f(x)$ when $|f(x)| \le \lambda$, $f_1(x) = 0$ when
$|f(x)| > 0$, $f_2(x) = f(x)$ when $|f(x)| > \lambda$, and $f_2(x) =
0$ when $|f(x)| \le \lambda$.  Notice in particular that $f_1$ lies in
$L^w$ for all $w \ge q$, and that $f_2$ lies in $L^u$ for all $u \le
q$.  For the present purposes, the idea is to use decompositions which
are better behaved, with $f_2$ having a more precise form as a sum of
multiples of atoms.  The Calder\'on--Zygmund method is again
applicable, although it should be mentioned that one first works
with $(p,r)$-atoms with one choice of $r$, and then afterwards
makes a conversion to a larger $r$ using the results described before.

	In addition to considering the effect of $T$ on atoms,
one can consider the effect of $T^*$ on atoms, and this leads
to conclusions about $T$ on $L^q$ for $q > 2$, by duality.


\end{document}